\documentclass[12pt]{article}
\usepackage{amsmath,amssymb,amsthm,booktabs}
\usepackage{mathrsfs}
\usepackage[OT1]{fontenc}
\usepackage[utf8]{inputenc}
\usepackage[colorlinks,citecolor=blue,urlcolor=blue]{hyperref}
\usepackage{txfonts}
\usepackage{indentfirst}
\usepackage{multirow}
\usepackage{float,subfig}

\usepackage{bm}
\usepackage{euscript}
\usepackage{graphicx}
\usepackage{multicol}
\usepackage[usenames,dvipsnames,svgnames,table]{xcolor}
\usepackage{lineno}

\numberwithin{equation}{section} \theoremstyle{plain}
\newtheorem{theorem}{Theorem}[section]
\newtheorem{lemma}{Lemma}[section]

\newtheorem{example}{Example}

\newtheorem{remark}{Remark}[section]

\def\blue{\color{blue}}
\numberwithin{equation}{section}

\begin{document}

\title{A varying terminal time mean-variance model
\footnotemark[1]}

\author{Shuzhen Yang
\footnotemark[2]
\footnotemark[3]}

\renewcommand{\thefootnote}{\fnsymbol{footnote}}

\footnotetext[1]{

\ \ \textbf{MSC2010 subject classification}: 91B28; 93E20; 49N10.

\ \ \textbf{OR/MS subject classification}: Finance/portfolio; dynamic programming/optimal control.}

\footnotetext[2]{Shandong University-Zhong Tai Securities Institute for Financial Studies, Shandong University, PR China, (yangsz@sdu.edu.cn).}
\footnotetext[3]{This work was supported by the National Key R\&D program of China (Grant No.2018YFA0703900), National Natural Science Foundation of China (Grant No.11701330), and Young Scholars Program of Shandong University.}

\date{}
\maketitle

 \textbf{Abstract}: To improve the efficient frontier of the classical mean-variance model in continuous time, we propose a varying terminal time mean-variance model with a  constraint on the mean value of the portfolio asset,  which moves  with the varying terminal time. Using  the embedding technique from stochastic optimal control in continuous time and by varying the terminal time, we determine an optimal strategy and related deterministic terminal time for the model. Our results suggest that doing so for  an  investment plan requires  minimizing  the variance with a varying terminal time.

\addcontentsline{toc}{section}{\hspace*{1.8em}Abstract}

\section{Introduction}

Since the seminal works of  {\blue Markowitz (1952, 1959)}, the mean-variance model has been used to balance the return (mean value) and risk (variance)  in a single-period portfolio selection model, (refer to
{\blue Markowitz (2014)}). Based on   mild assumptions, {\blue Merton (1972)} solved  this single-period problem analytically. More recently, multi-period and continuous time mean-variance portfolio selection models have been proposed. For example,  {\blue Richardson (1989)} studied a mean-variance model in a continuous-time setting for a single stock with a constant risk-free rate. {\blue Bajeux-Besnainou and Portait (1998)} considered dynamic asset allocation in a mean-variance framework. For the multi-period case,  {\blue Li and Ng
(2000)} embedded the discrete-time multi-period mean-variance problem within a multi-objective optimization framework. In the continuous time case,  {\blue Zhou and
Li (2000)} formulated the continuous-time mean-variance problem  as a stochastic linear-quadratic (LQ) optimal control problem. The solution to this problem is obtained by extending  the embedding technique introduced in  {\blue Li and Ng
(2000)} and using the results from the stochastic LQ control. Further extensions to the  mean-variance problem include those with bankruptcy prohibition, transaction costs, and random parameters in complete and incomplete markets ({\blue Bielecki, Jin, Pliska and Zhou (2005);
Dai, Xu and Zhou (2010); Lim (2004); Lim and Zhou (2002); Xia (2005)}).

For the aforementioned multi-period and continuous time cases, we derive the pre-committed strategies that differ from that of the single-period case; for further details, see {\blue   Kydland and Prescott (1997)}. {\blue Basak and Chabakauri (2010)} adopted a dynamic method to study  the mean-variance model and  {\blue Bj\"{o}rk, Murgoci and Zhou (2014)} studied the mean-variance problem with state dependent risk aversion.

In the classical mean-variance model, for a given deterministic terminal time $\tau$, we denote $X^{\pi}(\tau)$ as the terminal value of a portfolio asset with  strategy $\pi(\cdot)$, $\mathbb{E}[X^{\pi}(\tau)]$ and $\mathrm{Var} (X^{\pi}(\tau))=\mathbb{E}\big{(}X^{\pi}(\tau)-\mathbb{E}[X^{\pi}(\tau)]\big{)}^2$ as the mean and variance, respectively.
Note that we always want to minimize the variance $\mathrm{Var} (X^{\pi}(\tau))$ for a given mean level $\mathbb{E}[X^{\pi}(\tau)]=L$ in the single-period, multi-period, and continuous time cases, where $L$ is a constant. In the single-period case, $L$ can be viewed as the rate of return $\mathbb{E}[X^{\pi}(\tau)]$ over one period. However, in multi-period and continuous time cases, $L$ can be viewed only as  the return $\mathbb{E}[X^{\pi}(\tau)]$ over the terminal time $\tau$. We recognize this as an important difference between the single-period, multi-period case, and continuous time cases.


{  In the continuous time investment portfolio selection problem, we want to minimize the variance of $X^{\pi}(\cdot)$ for a given mean value of $X^{\pi}(\cdot)$ at some deterministic time $\tau$. However, an investor always needs to determine the length of the holding time before investing, for example, $\tau$ equals to $1$ or $2$ years. A natural question is weather we can choose the holding time $\tau>0$  as a strategy for the mean variance model in continuous time.  On the other hand, the investor would like to consider a short holding time $\tau$ within the same return. Note that the target of the mean value will move with the holding  time $\tau$ but not a constant. For any given deterministic  time $\tau>0$, we can stop the investment at a deterministic time $\tau^{\pi}\in[0,\tau]$, where $\tau^{\pi}$ depends on both the mean value $\mathbb{E}[X^{\pi}(\cdot)]$ and  time $\tau$. The criteria used to decide when to stop the investment is  as follows:
\begin{equation}
\label{intime-1}
\tau^{\pi}=\inf\bigg{\{}t:\mathbb{E}[X^{\pi}(t)]\geq xh(\tau),\ t\in [0,\tau] \bigg{\}},
\end{equation}
where $X^{\pi}(0)=x$ and $xh(\cdot)$ describes the target mean  of the asset $X^{\pi}(\cdot)$.
From the definition of $\tau^{\pi}$, we can see that $\mathbb{E}[X^{\pi}(t)]< xh(\tau),\ t<\tau^{\pi}$ and $\mathbb{E}[X^{\pi}(\tau^{\pi})]= xh(\tau)$. Thus, the criterion (\ref{intime-1}) shows that $\mathbb{E}[X^{\pi}(\cdot)]$ can take value $xh(\tau)$ at the smallest time $\tau^{\pi}$ for the given time $\tau>0$. }

Therefore, the investor  only need to minimize the variance at time $\tau^{\pi}$, and the cost functional is given as follows:
\begin{equation}
\label{incos-2}
J(\pi(\cdot))=\mathrm{Var} (X^{\pi}(\tau^{\pi}))=\mathbb{E}\big{(}X^{\pi}(\tau^{\pi})-\mathbb{E}[X^{\pi}(\tau^{\pi})]\big{)}^2.
\end{equation}
Here,   the terminal time $\tau^{\pi}$ depends on the control $\pi(\cdot)$ and time $\tau$, unlike  in the classical mean-variance problem. In this study, we will investigate an optimal strategy $\pi^*(\cdot)$  and a deterministic  terminal time $\tau^*$ for
the proposed model.

The work on the varying terminal time optimal control problem
that most closely related to ours is that of  {\blue Yang (2019)}, who established a stochastic maximum principle for a general state equation and cost functional. However, the result of {\blue Yang (2019)} is a necessary condition for an optimal strategy and, thus cannot be used to solve the varying terminal time mean-variance model in this study. { {\blue DeMiguel, Garlappi and Uppal (2009) } evaluated the out-of-sample performance of the sample-based mean-variance model, and compared with the naive 1/N portfolio. Based on the results of {\blue DeMiguel, Garlappi and Uppal (2009)}, we want to consider the effect of the number of risky assets to mean-variance model in this study. {\blue Dybvig (1988)} proposed the cost-efficient approach to the optimal portfolio selection in a straightforward way. Based on the cost-efficient approach, {\blue Bernard and  Vanduffel (2014)} considered the problem of mean-variance optimal portfolio in the presence of a benchmark, further see {\blue  Bernard, Vanduffel and Ye (2019)}. {\blue Bi, Jin and Meng (2018)} considered the behavioral mean–variance portfolio selection problem in continuous time. We point out that the results of {\blue Bernard and Vanduffel (2014)} can solve the first step of our problem which is consistent with the results of {\blue Zhou and Li (2000)}, see Example \ref{ex-1} in this paper.}

The remainder of this paper is organized as follows. In Section 2, we formulate the varying terminal time mean-variance model. Then, in Section 3, we investigate an optimal strategy and  the related terminal time for the proposed model. In Section 4, based on the main results of  Section 3, we compare our varying terminal time mean-variance model with the classical mean-variance model. Finally, we conclude  the  paper in Section 5.

\section{A new mean-variance model }

{ In the following, we will use $\tau$ to denote a given terminal time, $\tau^{\pi}$ to denote the terminal time related with strategy $\pi(\cdot)$, $\bar{\pi}(\cdot)$ to denote the optimal strategy under a given terminal  time,  and $(\pi^*(\cdot),\tau^*)$ to denote an optimal pair.}

Let $W$ be a $d$-dimensional standard Brownian motion defined in a complete
filtered probability space $(\Omega,\mathcal{F},P;\{ \mathcal{F}(t)\}_{t\geq
0})$, where $\{ \mathcal{F}(t)\}_{t\geq0}$ is the $P$-augmentation of the
natural filtration generated by  $W$. One risk-free  bond asset and $n$ risky  stock assets  are traded in the market. The bond satisfies the following equation:
\begin{eqnarray*}
\left\{\begin{array}{rl}
\mathrm{d}S_0(t) & \!\!\!= r(t)S_0(t)\mathrm{d}t,\;\;  t>0,\\
 S_0(0) & \!\!\!= s_0>0,
\end{array}\right.
\end{eqnarray*}
and the $i$'th ($1\leq i\leq n $) stock asset is described by
\begin{eqnarray*}
\left\{\begin{array}{rl} dS_i(t) & \!\!\!=
b_i(t)S_i(t)\mathrm{d}t+ \displaystyle S_i(t) \sum_{j=1}^{{d}}\sigma_{ij}(t)\mathrm{d}W_{j}(t),\;\;
t>0,\\
 S_i(0) & \!\!\!= s_i>0,\;\;
\end{array}\right.
\end{eqnarray*}
where $r(\cdot)\in \mathbb{R}$ is the risk-free return rate of the bond, $b(\cdot)=(b_1(\cdot),\cdots,b_n(\cdot))\in{\mathbb{R}^n}$ is the expected return rate of the risky asset, and $\sigma(\cdot)=(\sigma_{1}(\cdot),\cdots,\sigma_{n}(\cdot))^{\top}\in \mathbb{R}^{n\times d} $ is the corresponding volatility matrix. Given initial capital $x>0$, $\displaystyle \beta(\cdot)=(\beta_1(\cdot),\cdots,\beta_n(\cdot))\in \mathbb{R}^n$, where $\beta_i(\cdot)=b_i(\cdot)-r(\cdot),\ 1 \leq i \leq n$. The investor's wealth $X^{\pi}(\cdot)$ satisfies
\begin{equation}\label{asset}
\left\{\begin{array}{rl}
\!\!\!\mathrm{d}X^{\pi}(t)  & \!\!\!=\big[r(t)X^{\pi}(t)  +\beta(t)\pi(t)^{\top}   \big] \mathrm{d}t+\pi(t)\sigma(t)  \mathrm{d}W(t),  \\
 \!\!\!X^{\pi}(0) & \!\!\!=x,
\end{array}\right.
\end{equation}
where $\pi(\cdot)=(\pi_1(\cdot),\cdots,\pi_n(\cdot))\in \mathbb{R}^{n}$ is the capital invested in the risky asset $S(\cdot)=(S_1(\cdot),\cdots,S_n(\cdot))\in \mathbb{R}^n$ and  $\pi_0(\cdot)$ is the capital invested in the bond. Thus, we have
$
\displaystyle X^{\pi}(\cdot)=\sum_{i=0}^n\pi_i(\cdot).
$

In this study, we consider the following varying terminal time mean-variance model:
\begin{equation}
J(\pi(\cdot))=%
\mathbb{E}\big{(}X^{\pi}(\tau^{\pi})-\mathbb{E}[X^{\pi}(\tau^{\pi})]\big{)}^2,\label{cost-1}%
\end{equation}
with the following constraint on  the varying terminal time, for any given deterministic time $\tau>0$, we define the criteria for deciding when to stop the investment  as follows:
\begin{equation}
\label{time-1}
\tau^{\pi}=\inf\bigg{\{}t:\mathbb{E}[X^{\pi}(t)]\geq xh(\tau),\ t\in [0,\tau] \bigg{\}}.
\end{equation}
The definition of the deterministic terminal time $\tau^{\pi}$ shows that the mean value of wealth $X^{\pi}(\cdot)$ is moves with the target $xh(\cdot)$. Note that if $\tau^{\pi}<+\infty$, we have that $\mathbb{E}[X^{\pi}(\tau^{\pi})]= xh(\tau)$ and $ \mathbb{E}[X^{\pi}(t)]< xh(\tau),\ t<\tau^{\pi}$. This is the main difference between our model (\ref{cost-1}) with constraint (\ref{time-1})  and the classical mean-variance model.

The set of admissible strategies $\pi(\cdot)$ is defined as follows:
$$
\mathcal{A}=\bigg{\{}\pi(\cdot):  \pi(\cdot)\in L^2_{\mathcal{F}}[0,\tau^{\pi};\mathbb{R}^n],  \ \tau^{}\in (0,+\infty) \bigg{\}}.
$$
$L^2_{\mathcal{F}}[0,\tau;\mathbb{R}^n]$ is the set of all $\mathbb{R}^n$ valued, measurable processes $f(\cdot)$ adapted to $\{\mathcal{F}_t\}_{t\geq 0}$ such that
$$
\mathbb{E}\bigg{[}\int_0^{\tau}\left|f(t) \right|^2\mathrm{d}t \bigg{]}<+\infty.
$$
If there exists an optimal strategy $\pi^{*}(\cdot)
\in \mathcal{A}$ and the related deterministic terminal time $\tau^*$ that yields the minimum value of the cost functional (\ref{cost-1}), and if $\tau^*$ satisfies constraint  (\ref{time-1}), then we say $(\pi^{*}(\cdot),\tau^*)$ is the optimal strategy of the varying terminal time mean-variance model (\ref{cost-1}).

\begin{remark}\label{le-0}
Consider a special case of model (\ref{cost-1}) with constraint (\ref{time-1}). We suppose that the target of the mean value is a constant, $h(t)=L,\ t\geq 0$, where $L$ is a constant. Here, we want to determine   the optimal terminal time for a constant target. For any given $\tau>0$, applying the result of  {\blue Zhou and Li (2000)}, we  obtain the following representation for the variance,
\begin{equation*}
\displaystyle  \mathrm{Var}(X^{\bar{\pi}}(\tau))=
\frac{\bigg{(}\mathbb{E}[X^{\bar{\pi}}(\tau)]
-xe^{\int_0^{\tau}r(t)\mathrm{d}t}\bigg{)}^2}{e^{\int_0^{\tau}\phi(t)\mathrm{d}t}-1},
\end{equation*}
{ where $\phi(t)=\beta(t)[\sigma(t)\sigma(t)^{\top}]^{-1}\beta(t)^{\top}$}. Note that, $\mathbb{E}[X^{\bar{\pi}^{}}(\tau)]=L$. Thus, let
\begin{equation}\label{eq-0}
L
-xe^{\int_0^{\tau}r(t)\mathrm{d}t}=0.
\end{equation}
We further suppose that $\tau^*$  is the solution of equation (\ref{eq-0}), and that the related optimal strategy is $\pi^*(t)=0,\ t\leq \tau^*$.  Then, $\mathrm{Var}(X^{{\pi}^{*}}(\tau^*))=0$ implies that  we  invest all capital into the bound until time $\tau^*$. Thus, we need to provide the conditions for the function $h(\cdot)$ such that  $\mathbb{E}[X^{\bar{\pi}^{}}(\tau)]>xe^{\int_0^{\tau}r(s)\mathrm{d}s},\ \tau> 0$.
\end{remark}

We assume the following conditions,  which  we use  to obtain the optimal strategy for the proposed  model (\ref{cost-1}):

{$\textbf{H}_1$}: $r(\cdot), b(\cdot)$ and $\sigma(\cdot)$ are bounded deterministic continuous functions.

{$\textbf{H}_2$}: $r(\cdot),\beta(\cdot)>0$, $\sigma(\cdot)\sigma(\cdot)^{\top}>\varepsilon  \textbf{I}$.

{$\textbf{H}_3$}: $h(\cdot)$ is  increasing and differentiable in $[0,+\infty)$ with $h(0)> 1$, and satisfies
$$
\displaystyle {h(t)}> e^{\int_0^{t}r(s)\mathrm{d}s},\quad  t\geq 0,
$$
where    $\varepsilon>0$ is a given constant and $\textbf{I}$ is the identity matrix of $\mathbb{R}^{n\times n}$.

\begin{remark}\label{re-1} Note that $xh(\cdot)$ can be viewed  as the return of  asset $X^{{\pi}^{}}(\cdot)$, with initial wealth $x$.  By denoting $\displaystyle \psi(\cdot)= \frac{h'(\cdot)}{h(\cdot)}$, we  obtain
\[
{xh(t)}=  xh(0)e^{\int_0^{t}\psi(s)\mathrm{d}s}> xe^{\int_0^{t}r(s)\mathrm{d}s}, \quad t\geq0,
\]
yielding
$
\psi(t)>r(t)-\ln h(0).
$
Note that if $h(0)\leq 1$, we have $\tau^{\pi}=0$ for any strategy $\pi(\cdot)$.
\end{remark}

\section{Optimal strategy and  terminal time}

\subsection{Standard mean-variance model}
In this section, we investigate an optimal strategy $\pi(\cdot)$ for the problem defined in (\ref{cost-1}), with constraint (\ref{time-1}) for varying terminal time $\tau^{\pi}$.
Here,  we  describe how to construct an optimal strategy for (\ref{cost-1}) with varying terminal time (\ref{time-1}). The  steps are as follows:

\textbf{Step 1}: For any given deterministic  $\tau\in (0,+\infty)$, we determine an optimal strategy $\bar{\pi}^{}(\cdot)$ for the related mean-variance model using the embedding technique proposed  in {\blue Zhou and Li (2000)}.

\textbf{Step 2}: Verify that the optimal strategy $\bar{\pi}^{}(\cdot)$  from \textbf{Step 1} is an element of the set of admissible strategies  $\mathcal{A}$ and ${\tau}^{\bar{\pi}^{}}=\tau$. Thus $\mathcal{A}\neq \varnothing$.

\textbf{Step 3}: Minimize the variance over $\tau\in (0,+\infty)$ to obtain, the optimal strategy $\pi^*(\cdot)$ and the related terminal time $\tau^{*}$ for problem (\ref{cost-1}) with varying terminal time (\ref{time-1}).

\bigskip

Next, we consider  \textbf{Step 1} in greater detail. For any given deterministic time $\tau\in(0,+\infty)$, we introduce the following mean-variance problem: minimize the  cost functional,
\begin{equation}
\label{cos-step-1}
J(\pi^{}(\cdot))=-\mathbb{E}[X^{\pi^{}}(\tau)]+\frac{\mu(\tau)}{2} \mathrm{Var} (X^{\pi^{}}(\tau)).
\end{equation}
To solve the cost functional (\ref{cos-step-1}), we employ the following model:
 \begin{equation}
\label{cos-step-2}
J(\pi^{}(\cdot))=\mathbb{E}[\frac{\mu(\tau)}{2} X^{\pi^{}}(\tau)^2-\lambda(\tau) X^{\pi^{}}(\tau)].
\end{equation}
Applying the embedding technique of {\blue Zhou and Li (2000)} for mean-variance models in  the continuous time case, we have the following results. For  $  t\leq \tau$, and denoting
\begin{equation*}
\begin{array}
[c]{rl}%
& \phi(t)=\beta(t)[\sigma(t)\sigma(t)^{\top}]^{-1}\beta(t)^{\top};\\
&\bar{\lambda}(\tau)=e^{\int_0^{\tau}\phi(t)\mathrm{d}t}+\mu(\tau) x e^{\int_0^{\tau}r(t)\mathrm{d}t};\\
& \displaystyle \gamma(\tau)=\frac{\bar{\lambda}(\tau)}{\mu(\tau)}=
\frac{1}{\mu(\tau)}e^{\int_0^{\tau}\phi(t)\mathrm{d}t}+x e^{\int_0^{\tau}r(t)\mathrm{d}t},
\end{array}
\end{equation*}
the optimal strategy for (\ref{cos-step-1}) is given as follows:
$$
\bar{\pi}^{}(t)=[\sigma(t)\sigma(t)^{\top}]^{-1}\beta(t)^{\top}[\gamma(\tau) e^{-\int_t^{\tau}r(s)\mathrm{d}s}-X^{\bar{\pi}^{}}(t)],\quad  t\leq \tau.
$$
Here, $\mathbb{E}[{X}^{\bar{\pi}^{}}(\cdot)]$ and $\mathbb{E}[{X}^{\bar{\pi}^{}}(\cdot)^2]$ satisfy the following linear ordinary differential equations:
\begin{equation}\label{step-asset-1}
\left\{\begin{array}{rl}
\!\!\!\mathrm{d}\mathbb{E}[X^{\bar{\pi}^{}}(t)]  & \!\!\!=\bigg[r(t)\mathbb{E}[{X}^{\bar{\pi}^{}}(t)]  +\gamma (\tau) e^{-\int_t^{\tau}r(s)\mathrm{d}s}\phi(t)  \bigg] \mathrm{d}t,  \\
 \!\!\!\mathbb{E}[X^{\bar{\pi}^{}}(0)] & \!\!\!=x,
\end{array}\right.
\end{equation}
and
\begin{equation}\label{step-asset-2}
\left\{\begin{array}{rl}
\!\!\!\mathrm{d}\mathbb{E}[X^{\bar{\pi}^{}}(t)^2]  & \!\!\!=\bigg[(2r(t)-\phi(t))\mathbb{E}[{X}^{\bar{\pi}^{}}(t)^2]  +\gamma (\tau)^2 e^{-\int_t^{\tau}2r(s)\mathrm{d}s}\phi(t)  \bigg] \mathrm{d}t,  \\
 \!\!\!\mathbb{E}[X^{\bar{\pi}^{}}(0)^2] & \!\!\!=x^2.
\end{array}\right.
\end{equation}
Note that $\displaystyle \gamma(\tau)=
\frac{1}{\mu(\tau)}e^{\int_0^{\tau}\phi(t)\mathrm{d}t}+x e^{\int_0^{\tau}r(t)\mathrm{d}t}$; thus, we have
$$
\displaystyle \mathbb{E}[{X}^{\bar{\pi}^{}}(t)]=xe^{\int_0^t(r(s)-\phi(s))
\mathrm{d}s}+\gamma(\tau)e^{-\int_t^{\tau}r(s)
\mathrm{d}s}[1-e^{-\int_0^{t}\phi(s)\mathrm{d}s}].
$$
Let $t=\tau$. Then,  it follows that,
$$
\displaystyle \mathbb{E}[{X}^{\bar{\pi}^{}}(\tau)]=
\frac{1}{\mu(\tau)}(e^{\int_0^{\tau}\phi(t)\mathrm{d}t}-1)+x e^{\int_0^{\tau}r(t)\mathrm{d}t}.
$$
By assumption $\textbf{H}_3$,  we have $xh(\tau)-x e^{\int_0^{\tau}r(t)\mathrm{d}t}>0$. Thus, we can choose $\mu(\tau)>0$ such that
$
\mathbb{E}[{X}^{\bar{\pi}^{}}(\tau)]=xh(\tau),
$
where $\mu(\tau)$ satisfies
\begin{equation}\label{tau-1}
\displaystyle \mu(\tau)=\frac{e^{\int_0^{\tau}\phi(t)\mathrm{d}t}-1}{xh(\tau)-x e^{\int_0^{\tau}r(t)\mathrm{d}t}}.
\end{equation}

Based on the explicit solution for $\mathbb{E}[{X}^{\bar{\pi}^{}}(\cdot)]$, the following theorem presents the main results from  \textbf{Step 1}.
\begin{theorem}\label{the-1}
Let Assumptions $\textbf{H}_1, \textbf{H}_2,$ and $\textbf{H}_3$ hold. For any given deterministic time $\tau>0$, $(\bar{\pi}^{}(\cdot),X^{\bar{\pi}^{}}(\cdot))$ is an optimal pair for the mean-variance problem in (\ref{cos-step-1}). The efficient frontier is given as follows:
\begin{equation}\label{step-the}
\displaystyle  \mathrm{Var}(X^{\bar{\pi}^{}}(\tau))=
\frac{\bigg{(}\mathbb{E}[X^{\bar{\pi}^{}}(\tau)]
-xe^{\int_0^{\tau}r(t)\mathrm{d}t}\bigg{)}^2}{e^{\int_0^{\tau}\phi(t)\mathrm{d}t}-1}.
\end{equation}
\end{theorem}

{\begin{example}\label{ex-1}
 In \textbf{Step 1}, we use the embedding technique of   {\blue Zhou and Li (2000)} to solve the mean-variance model in continuous time. Now, we consider a one-dimensional Black-Scholes setting: there is bond with a constant risk-free rate $r>0$, and a risky stock $S$:
\begin{eqnarray*}
\left\{\begin{array}{rl} dS(t) & \!\!\!=
\mu S(t)\mathrm{d}t+ \sigma\displaystyle S(t) \mathrm{d}W(t),\;\;
t>0,\\
 S(0) & \!\!\!= 1>0.\;\;
\end{array}\right.
\end{eqnarray*}
By Theorem \ref{the-1}, we have
\begin{equation}\label{ex-eq-1}
\displaystyle  \mathrm{Var}(X^{\bar{\pi}^{}}(\tau))=
\frac{\bigg{(}\mathbb{E}[X^{\bar{\pi}^{}}(\tau)]
-xe^{r \tau }\bigg{)}^2}{e^{\phi^2 \tau }-1},
\end{equation}
where $\displaystyle \phi=\frac{\mu-r}{\sigma}$.

Applying Proposition 3.1 of {\blue Bernard and Vanduffel (2014)}, also see Subsection 3.2 in {\blue Bernard and Vanduffel (2014)},  we have
\begin{equation}\label{ex-eq-2}
X^{\bar{\pi}^{}}(\tau)=a-\frac{b}{S^*_{\tau}},
\end{equation}
where
$$
a=xe^{r\tau}+e^{\phi^2\tau}\sqrt{\frac{\mathrm{Var}(X^{\bar{\pi}^{}}(\tau))}{e^{\phi^2\tau}-1}},\
b=e^{r\tau}\sqrt{\frac{\mathrm{Var}(X^{\bar{\pi}^{}}(\tau))}{e^{\phi^2\tau}-1}},
$$
and
$$
S_{\tau}^*=e^{\big(\frac{\phi \mu}{\sigma}+(1-\frac{\phi}{\sigma})r-\frac{\phi^2}{2}\big)\tau+\phi W(\tau)}.
$$
It follows that $\mathbb{E}[\frac{1}{S_{\tau}^*}]=e^{-r\tau}$. With expectation on both sides of equation (\ref{ex-eq-2}), and by plugging $\mathbb{E}[\frac{1}{S_{\tau}^*}],a,b$ into it, one obtains,
$$
\big(\mathbb{E}[X^{\bar{\pi}^{}}(\tau)]-xe^{r\tau}\big)^2={\mathrm{Var}(X^{\bar{\pi}^{}}(\tau))}({e^{\phi^2\tau}-1}),
$$
which is same with equation (\ref{ex-eq-1}) and  is consistent with the results of Theorem \ref{the-1}.

{\blue Bernard and Vanduffel (2014)} applied the cost-efficient
approach to optimal portfolio selection in a straightforward way, and obtained the optimal payoff $X^{\bar{\pi}^{}}(\tau)$. Furthermore, we have
\begin{equation}\label{ex-eq-3}
X^{\bar{\pi}^{}}(\tau)=xe^{r\tau}+\big(e^{\phi^2\tau}-e^{-\phi W(\tau)-\frac{\phi^2\tau}{2}}\big)\sqrt{\frac{\mathrm{Var}(X^{\bar{\pi}^{}}(\tau))}{e^{\phi^2\tau}-1}}.
\end{equation}
From equation (\ref{ex-eq-3}), the optimal payoff can be divided into two parts: risk-free return
$xe^{r\tau}$, risky return $\big(e^{\phi^2\tau}-e^{-\phi W(\tau)-\frac{\phi^2\tau}{2}}\big)\sqrt{\frac{\mathrm{Var}(X^{\bar{\pi}^{}}(\tau))}{e^{\phi^2\tau}-1}}$, which shows that high return comes with high risk.

\end{example}
}

Next, we consider  \textbf{Step 2} in further details, and show that $\mathcal{A}\neq \varnothing$.
\begin{lemma}\label{le-3}
Let Assumptions $\textbf{H}_1, \textbf{H}_2$, and $\textbf{H}_3$ hold. Then, we have $\bar{\pi}^{}(\cdot)\in \mathcal{A}$, ${\tau}^{\bar{\pi}^{}}=\tau$, and $\mathcal{A}\neq \varnothing$.
\end{lemma}
\noindent \textbf{Proof:} The proof is given in Appendix A. $\ \ \ \ \ \ \ \ \Box$

\bigskip

Lastly, we examine \textbf{Step 3} more closely. We want to obtain an optimal strategy $\pi^*(\cdot)$ for  model (\ref{cost-1}) and the related deterministic terminal time $\tau^*$. For simplicity of notation,  we consider the following decomposition of $xh(\cdot)$,
\begin{equation}\label{eqh-0}
xh(\tau)=x(h(0)-1) e^{\int_0^{\tau}\frac{\theta(t)}{2}\mathrm{d}t}+ xe^{\int_0^{\tau}r(t)\mathrm{d}t},\ \tau> 0,
\end{equation}
where
$$\displaystyle \theta(s)=\frac{2h'(s)-2r(s)e^{\int_0^{s}r(t)\mathrm{d}t}}
{h(s)-e^{\int_0^{s}r(t)\mathrm{d}t}},\ s\geq 0.$$
Here, $x(h(0)-1) e^{\int_0^{\tau}\frac{\theta(t)}{2}\mathrm{d}t}$ can be viewed as the excess return of $\mathbb{E}[X^{\bar{\pi}^{}}(\tau)]$ over the return of bond $xe^{\int_0^{\tau}r(t)\mathrm{d}t}$, where $\bar{\pi}(\cdot)$ in $\mathbb{E}[X^{\bar{\pi}^{}}(\tau)]$ denotes the optimal strategy with terminal time $\tau$.
\begin{remark}\label{re0-2}
For any given $\tau_1,\tau_2> 0$, we have
\begin{equation*}
\begin{array}
[c]{rl}%
\displaystyle &\mathbb{E}[{X}^{\bar{\pi}^{}}(\tau_1)]-xe^{\int_0^{\tau_1}r(t)\mathrm{d}t}\\
=&x(h(0)-1) e^{\int_0^{\tau_1}\frac{\theta(t)}{2}\mathrm{d}t}\\
=&x(h(0)-1)e^{\int_0^{\tau_2}\frac{\theta(t)}{2}\mathrm{d}t}e^{-\int_{\tau_1}^{\tau_2}\frac{\theta(t)}{2}\mathrm{d}t}\\
=&\bigg{(}\mathbb{E}[{X}^{\bar{\pi}^{}}(\tau_2)]-xe^{\int_0^{\tau_2}r(t)\mathrm{d}t}\bigg{)}e^{-\int_{\tau_1}^{\tau_2}\frac{\theta(t)}{2}\mathrm{d}t}.
\end{array}
\end{equation*}
Based on the discounted factor $e^{-\int_{\tau_1}^{\tau_2}\frac{\theta(t)}{2}\mathrm{d}t}$,   we call that $\mathbb{E}[{X}^{\bar{\pi}^{}}(\tau_1)]$ and $\mathbb{E}[{X}^{\bar{\pi}^{}}(\tau_2)]$ have the same excess return. Therefore,  we can compare the variance of the asset ${X}^{\bar{\pi}^{}}(\tau)$ with different terminal time $\tau>0$ under the same excess return. In the following,  we show the existence of the minimum value of $\mathrm{Var}(X^{\bar{\pi}^{}}(\tau))$ over $\tau\in (0,+\infty)$ under the  same excess return.
\end{remark}

\begin{lemma}\label{le-5}
Let Assumptions $\textbf{H}_1, \textbf{H}_2$, and $\textbf{H}_3$ hold. We have the following results:

(i). { If
$$
\theta(s)> \phi(s)+\delta ,\ s\geq 0,
$$
for a constant $\delta>0$ and }
$$\displaystyle \theta(s)=\frac{2h'(s)-2r(s)e^{\int_0^{s}r(t)\mathrm{d}t}}
{h(s)-e^{\int_0^{s}r(t)\mathrm{d}t}},\ s\geq 0,$$
there exist $\pi^*(\cdot)$ and related deterministic terminal time $\tau^*\in (0,+\infty)$ such that
\begin{equation*}
\displaystyle  \mathrm{Var}(X^{{\pi}^{*}}(\tau^*))=\min_{\tau\in(0,+\infty)}\mathrm{Var}(X^{\bar{\pi}^{}}(\tau)),
\end{equation*}
where $\bar{\pi}(\cdot)$ in $\mathbb{E}[X^{\bar{\pi}^{}}(\tau)]$ denotes the optimal strategy with terminal time $\tau$, $\mathbb{E}[{X}^{{\pi}^{*}}(\tau^*)]=xh(\tau^*)$ and
\begin{equation*}
\begin{array}
[c]{rl}%
& \phi(t)=\beta(t)[\sigma(t)\sigma(t)^{\top}]^{-1}\beta(t)^{\top},\quad  t\leq \tau^*;\\
& \displaystyle \gamma(\tau^*)=
\frac{1}{\mu(\tau^*)}e^{\int_0^{\tau^*}\phi(t)\mathrm{d}t}+x e^{\int_0^{\tau^*}r(t)\mathrm{d}t};\\
& \pi^*(t)=[\sigma(t)\sigma(t)^{\top}]^{-1}\beta(t)^{\top}(\gamma(\tau^*) e^{-\int_t^{\tau^*}r(t)\mathrm{d}t}-X^{{\pi}^{*}}(t)),\quad  t\leq \tau^*.
\end{array}
\end{equation*}

(ii).  If
$$
\theta(s)\leq \phi(s) ,\ s\geq  0,
$$
we have that
$$
\lim_{\tau\to +\infty}\mathrm{Var}(X^{\bar{\pi}^{}}(\tau))=\inf_{\tau\in(0,+\infty)}
\mathrm{Var}(X^{\bar{\pi}^{}}(\tau)),
$$
where
$$
\bar{\pi}^{}(t)=[\sigma(t)\sigma(t)^{\top}]^{-1}\beta(t)^{\top}(\gamma(\tau) e^{-\int_t^{\tau}r(t)\mathrm{d}t}-X^{\bar{\pi}^{}}(t)),\quad  t\leq \tau.
$$
\end{lemma}

\noindent \textbf{Proof:} The proof is given in Appendix A. $\ \ \ \ \ \ \ \ \Box$

\begin{remark}
The results of Lemma \ref{le-5} show  that if we consider the mean-variance model for a given terminal time, the variance may not determine the minimum value. Thus, we consider a varying terminal time mean-variance model with a constraint on the varying terminal time.  In the following, we prove that the optimal strategy $\pi^*(\cdot)$ and optimal terminal time $\tau^*$ in Lemma  \ref{le-5} solve model (\ref{cost-1}) with constraint (\ref{time-1}).
\end{remark}

\begin{theorem}\label{the-2}
Let Assumptions $\textbf{H}_1, \textbf{H}_2$, and $\textbf{H}_3$ hold. Then, we have the following results:

(i). { If
$$
\theta(s)> \phi(s)+\delta ,\ s\geq 0,
$$
for a constant $\delta>0$}, then $(\pi^*(\cdot),\tau^*)$  given in Lemma \ref{le-5} solves the varying terminal time mean-variance model (\ref{cost-1}) with constraint (\ref{time-1}).

(ii).  If
$$
\theta(s)\leq \phi(s) ,\ s\geq 0,
$$
for the optimal strategy $(\bar{\pi}^{}(\cdot),\tau)$ of cost functional (\ref{cos-step-1}), (\ref{cost-1}) is decreasing with $\tau\in (0,+\infty)$. Thus, the proposed   model does not
  yield an optimal strategy in finite time.
\end{theorem}

\noindent \textbf{Proof:} The proof is given in Appendix A. $\ \ \ \ \ \ \ \ \Box$

\begin{remark}\label{re-4}
Comparing with results of Lemma \ref{le-5} and Theorem \ref{the-2},  we find that the optimal strategy $(\pi^*(\cdot),\tau^*)$ of  $\displaystyle \min_{\tau\in(0,+\infty)}\mathrm{Var}(X^{\bar{\pi}^{}}(\tau))$   is the same with the optimal strategy of model (\ref{cost-1}) with  constraint (\ref{time-1}). For any given $\tau> 0$, $\mathrm{Var}(X^{{\pi}}(\tau))$ takes the minimum value at strategy $(\bar{\pi}^{}(\cdot),\tau)$, with $\mathbb{E}[X^{\bar{\pi}^{}}(\tau)]-xh(\tau)=0$.  Note that in the constraint
$$
\tau^{\pi}=\inf\bigg{\{}t:\mathbb{E}[X^{\pi}(t)]\geq xh(\tau),\ t\in [0,\tau] \bigg{\}},
$$
which includes all  strategies $\pi(\cdot)$ such that $\mathbb{E}[X^{{\pi}}(t)]-xh(t)=0$
at $t=\tau^{\pi}>0$, where $\tau^{\pi}$ is the minimum time such that $\mathbb{E}[X^{{\pi}}(\tau^{\pi})]-xh(\tau^{\pi})=0$. Lemma \ref{le-3} shows that $(\bar{\pi}^{}(\cdot),\tau)$ satisfies  constraint (\ref{time-1}). We conclude that the optimal strategy of the classical mean-variance model  is consistent with
constraint (\ref{time-1}) at a given terminal time $\tau> 0$.
\end{remark}

Based on Theorem \ref{the-2}, we compare the varying terminal model mean-variance model  (\ref{cost-1}) with classical mean-variance model (\ref{cos-step-1}) under a given deterministic terminal time $\tau>0$ in the following section.

\subsection{Multi-dimensional Black-Scholes setting}

{
Now, we apply the results of Theorem \ref{the-2} to consider a multi-dimensional Black-Scholes market. The bond satisfies the following equation:
\begin{eqnarray*}
\left\{\begin{array}{rl}
\mathrm{d}S_0(t) & \!\!\!= rS_0(t)\mathrm{d}t,\;\;  t>0,\\
 S_0(0) & \!\!\!= s_0>0,
\end{array}\right.
\end{eqnarray*}
and the price of the $i$'th ($1\leq i\leq n$) stock asset is given by
\begin{eqnarray*}
\left\{\begin{array}{rl} dS_i(t) & \!\!\!=
b_iS_i(t)\mathrm{d}t+ \sigma_i\displaystyle S_i(t) \mathrm{d}W_i(t),\;\;
t>0,\\
 S_i(0) & \!\!\!= s_i>0,\;\;
\end{array}\right.
\end{eqnarray*}
 where $r,b_i,\sigma_i\in \mathbb{R}, 1\leq i \leq n$ are deterministic  parameters and independent from time $t>0$, and $(W_1(\cdot),W_2(\cdot),\cdots,W_n(\cdot))$ is a given $n$-dimensional Brownian motion.

Let $h(t)=\alpha e^{\int_0^t\frac{{\theta(s)}}{2}\mathrm{d}s}+e^{ rt}, \ t\geq 0$, where $\alpha$ are  given constants. In the following, we show the assumption for the parameters $r,\alpha,b_i,\sigma_i, 1\leq i \leq n$ and function $\theta(\cdot)$ as follows:

\bigskip

{$\textbf{H}_4$}: $r,\sigma_i,\alpha>0$, $b_i>r,\ 1\leq i\leq n$, $\theta(t)>0,\ t\geq 0$.

\begin{remark}\label{re-5} Notice that in Assumption $\textbf{H}_4$, we suppose $r,\sigma_i>0$, $b_i>r,\ 1\leq i\leq n$  which is natural in the financial market. Here, $\alpha>0$ means that the investor wants to obtain excess return  exceeding the return of the asset bond. The condition on $\theta(\cdot)$ shows that   the target of the mean value is increasing with the holding time.
\end{remark}
In the following, we assume Assumption $\textbf{H}_4$ holds in the remainder of this paper. In addition, if Assumption $\textbf{H}_4$ holds, we can show that Assumptions $\textbf{H}_1,\textbf{H}_2,\textbf{H}_3$ are right. In the representation of $xh(t)=x\alpha e^{\int_0^t\frac{\theta(s)}{2}\mathrm{d}s}+xe^{ rt},\ t\geq 0$,  the term $x\alpha e^{\int_0^t\frac{{\theta(s)}}{2}\mathrm{d}s}$  can be viewed as the excess return exceeding the return of asset bond $xe^{ rt}$. For any given deterministic time $\tau> 0$, and denoting
\begin{equation*}
\begin{array}
[c]{rl}%
&\displaystyle  \phi=\sum_{i=1}^n\bigg{(}\frac{b_i-r}{\sigma_i}\bigg{)}^2,\quad \displaystyle \gamma(\tau)=
\frac{1}{\mu(\tau)}e^{\phi\tau}+x e^{r\tau}, \quad
\displaystyle \displaystyle \mu(\tau)=\frac{e^{\phi \tau}-1}{x\alpha e^{\int_0^{\tau}\frac{\theta(s)}{2}\mathrm{d}s}}.
\end{array}
\end{equation*}
\begin{remark}\label{re-6}
Notice that, $xh(\cdot)$ is the objective of mean value $\mathbb{E}[{X}^{{\pi}^{}}(\cdot)]$. The term $xe^{ rt}$ denotes the return from bond at time $t$, while
$x\alpha e^{\int_0^t\frac{\theta(s)}{2}\mathrm{d}s}$ represents the excess return exceeding the return of asset bond at time $t$. In addition, we have
$$
\displaystyle \frac{\theta(t)}{2}=\frac{\big(x\alpha e^{\int_0^t\frac{\theta(s)}{2}\mathrm{d}s}\big)'}{x\alpha e^{\int_0^t\frac{\theta(s)}{2}\mathrm{d}s}}
$$
 is the rate of excess return at time $t$. $\displaystyle \phi=\sum_{i=1}^n\bigg{(}\frac{b_i-r}{\sigma_i}\bigg{)}^2$ is increasing with $n$ which is the number of risky assets, where $\displaystyle \frac{b_i-r}{\sigma_i}$ is the sharp ratio of the $i$'th risky asset, $1\leq i\leq n$.
\end{remark}

Therefore, the optimal strategy is
$$
\displaystyle \bar{\pi}^{}(t)=\big(\frac{b_1-r}{\sigma_1^2},\frac{b_2-r}{\sigma_2^2} ,\cdots,\frac{b_n-r}{\sigma_n^2} \big)^{\top} (\gamma(\tau)e^{-r\tau}-X^{\bar{\pi}^{}}(t)),\quad t\leq \tau,
$$
where $\bar{\pi}(\cdot)$  denotes the optimal strategy with terminal time $\tau$, and
$
\displaystyle \mathbb{E}[{X}^{\bar{\pi}^{}}(\tau)]=
\frac{e^{\phi \tau}-1}{\mu(\tau)}+x e^{r\tau}.
$
We can choose $\mu(\tau)>0$ such that
$
\displaystyle \mathbb{E}[{X}^{\bar{\pi}^{}}(\tau)]=x\alpha e^{\int_0^{\tau}\frac{\theta(s)}{2}\mathrm{d}s}+xe^{ r\tau}.
$

For the given  $\tau> 0$, the efficient frontier for the classical mean-variance model (\ref{cos-step-1}) is given as follows:
\begin{equation}\label{sa-ef}
\displaystyle  \mathrm{Var}(X^{\bar{\pi}^{}}(\tau))=
\frac{\bigg{(}\mathbb{E}[X^{\bar{\pi}^{}}(\tau)]
-xe^{r\tau}\bigg{)}^2}{e^{\phi \tau}-1}=\frac{x^2\alpha^2 e^{\int_0^{\tau}{\theta(s)}\mathrm{d}s}}{e^{\phi \tau}-1},
\end{equation}
where
$$
\displaystyle \mathbb{E}[{X}^{\bar{\pi}^{}}(\tau)]=x\alpha e^{\int_0^{\tau}\frac{\theta(s)}{2}\mathrm{d}s}+xe^{ r\tau}.
$$
Notice that, Remark \ref{re-6} shows that $\phi$ is increasing with $n$, which indicates that the variance $\mathrm{Var}(X^{\bar{\pi}^{}}(\tau))$ is decreasing with the number of risky assets. We consider the relation of $\mathrm{Var}(X^{\bar{\pi}^{}}(\tau))$ and $\displaystyle\frac{\theta(\tau)}{2}$ as the efficient frontier.

If $\displaystyle {\theta}(t)\leq  \phi,\ t\geq 0$, Lemma \ref{le-5} indicates that we cannot obtain an optimal terminal time. Thus, in the following, we suppose
$\displaystyle {\theta}(t)> \phi+\delta,\ t\geq 0$, for a constant $\delta>0$.
Based on the proof of Lemma \ref{le-5}, we have that the optimal terminal time is  the smallest solution of the following equation:
\begin{equation}\label{eq-1}
({\theta}(\tau)-\phi)e^{\phi\tau} -{\theta(\tau)}=0.
\end{equation}
\begin{remark}\label{re-7}
Notice that, $\displaystyle \frac{\theta(\cdot)}{2}$ is  the rate of excess return of the objective of mean value $\mathbb{E}[{X}^{{\pi}^{}}(\cdot)]$. Thus, the investor needs to give the value of $\theta(\cdot)$ before investing. In general, it is convenient for the investor  to suppose a constant rate of excess return. On the other hand,  if $\theta(\cdot)$ is a constant, it is easy to check the condition in Lemma \ref{le-5} and we can obtain an explicit solution for the optimal holding time $\tau^*$.
\end{remark}
}

 In the following, we suppose that ${\theta}(\tau)=\theta$ is a constant that is independent from $\tau$. We can obtain a unique solution for equation (\ref{eq-1}), as follows:
$$
\tau^*=\displaystyle \frac{1}{\phi}\ln\frac{{\theta}}{{\theta}-\phi},\quad
\phi=\sum_{i=1}^n\bigg{(}\frac{b_i-r}{\sigma_i}\bigg{)}^2.
$$
Therefore, the efficient frontier for  model (\ref{cost-1}) is given as follows:
\begin{equation}\label{sa-ef}
\displaystyle  \mathrm{Var}(X^{{\pi}^*}(\tau^*))=\frac{x^2\alpha^2 e^{{\theta}\tau^*}}{e^{\phi \tau^*}-1}=x^2\alpha^2\kappa\bigg{(}
\frac{\kappa}{\kappa-1}\bigg{)}^{\kappa-1},
\end{equation}
where
$$
\displaystyle \kappa=\frac{{\theta}}{\phi},\quad \mathbb{E}[X^{{\pi^*}}(\tau^*)]=x\alpha e^{\frac{{\theta}}{2}\tau^*}+xe^{ r\tau^*}.
$$
 Here, we consider the relation of $\mathrm{Var}(X^{\pi^{*}}(\tau^*))$ and $\displaystyle\frac{\theta}{2}$ as the efficient frontier.

\section{Compare with classical mean-variance model}

Following the results in Section 3, we  compare the variance $ \mathrm{Var}(X^{{\pi}^*}(\tau^*))$ with $\mathrm{Var}(X^{\bar{\pi}^{}}(\tau))$.
Based on Assumption $\textbf{H}_4$, we first consider $n=1$,  and take
\begin{equation*}
\begin{array}
[c]{rl}%
&\displaystyle  r=5\%,\ b_1=10\%,\ \sigma_1=20\%,\ \alpha_1=0.5,\ \alpha_2=0.3,\  {\theta}=40\%;\\
&y_1(\tau)=\displaystyle  \mathrm{Var}(X^{\bar{\pi}^{}}(\tau))=\frac{x^2\alpha_1^2 e^{{\theta}\tau}}{e^{\phi \tau}-1};\\
& \displaystyle y_2(\tau)=\displaystyle  \mathrm{Var}(X^{\bar{\pi}^{}}(\tau))=\frac{x^2\alpha_2^2 e^{{\theta}\tau}}{e^{\phi \tau}-1}.
\end{array}
\end{equation*}
In addition, we have that
$$
\displaystyle \kappa=\frac{{\theta}}{\phi}=6.4 >1,
$$
which satisfies  the case $(i)$ of Theorem \ref{the-2}. Notice that, $y_1(\tau)$
and $y_2(\tau)$ are the variances of the classical mean-variance model and  the functions of a given terminal time $\tau$ with $\alpha_1=0.5,\ \alpha_2=0.3$, respectively.

\begin{figure}[H]
\caption{ $\mathrm{Var}(X^{\bar{\pi}^{}}(\tau))$ varies with the terminal time $\tau>0$}
 \label{fig1}
\begin{center}
\includegraphics[width=4.2 in]{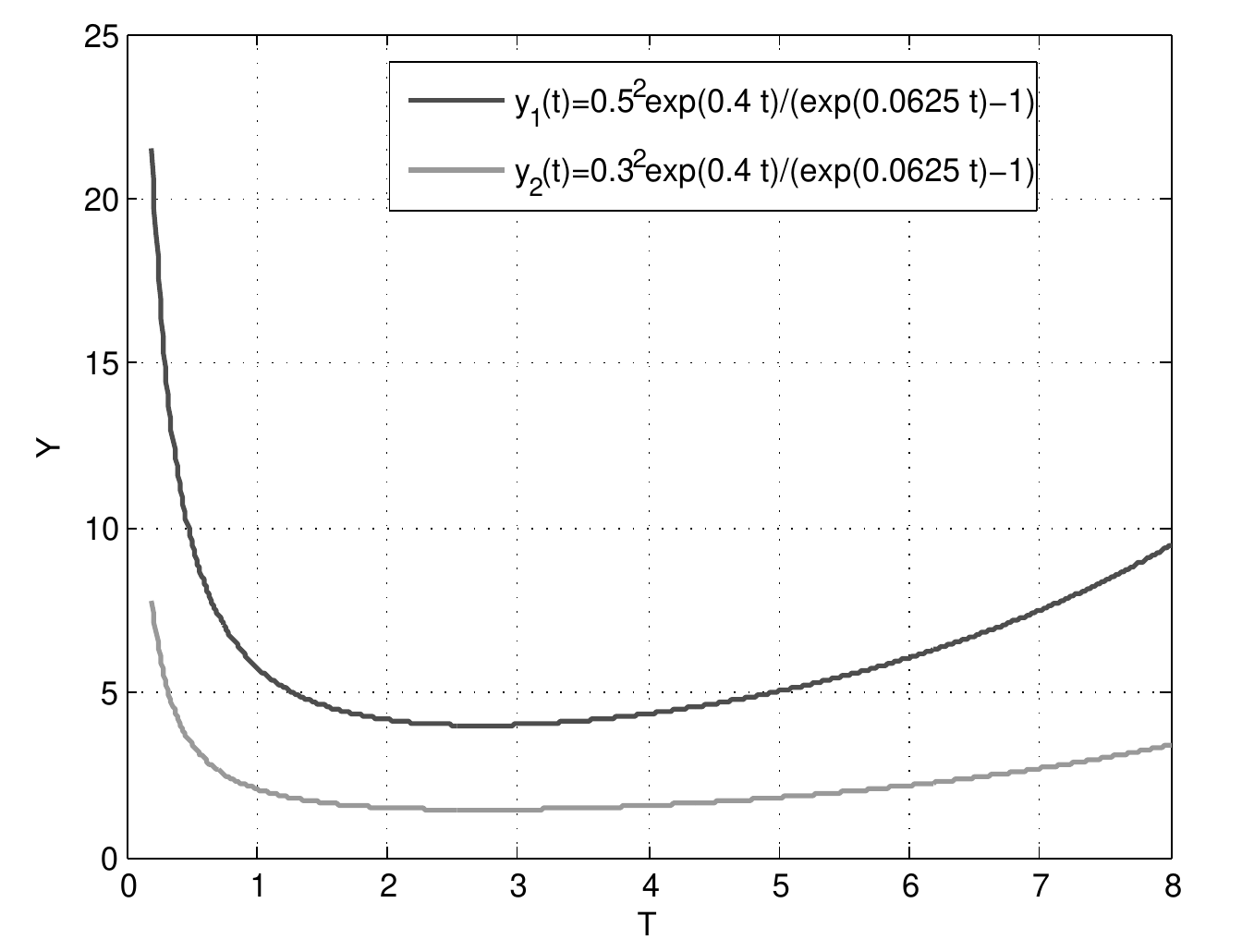}
\end{center}
\end{figure}
 In Figure \ref{fig1}, we plot the functions of $y_1(\cdot)$ and $y_2(\cdot)$ over $t\in [0.1,8]$ . We can verify that  $y_1(\cdot)$ and $y_2(\cdot)$ take  minimum values at $\tau^*=2.72$, and $\tau^*$ is given by
$
\tau^*=\displaystyle \frac{1}{\phi}\ln\frac{{\theta}}{{\theta}-\phi}.
$
For the given parameters $r=5\%,\ b_1=10\%,\ \sigma_1=20\%,\  {\theta}=40\%$, from Figure \ref{fig1},  we can see that $\mathrm{Var}(X^{\bar{\pi}^{}}(\tau))$ is decreasing in $[0.1,\tau^*]$ and increasing in $[\tau^*,8]$. Note that $h(0)=\alpha+1>1$ and  the optimal strategy $\bar{\pi}^{}(\cdot)$ yields,  a high risk (variance) for a given small terminal time $\tau<\tau^*$. At the same time, the optimal strategy $\bar{\pi}^{}(\cdot)$ of the investor could  yield a high risk (variance) for a given big terminal time $\tau>\tau^*$. These results show that we need to choose an optimal terminal time $\tau^*$ for a given rate of excess return $\displaystyle \frac{\theta}{2}$. In addition, we show that $\mathrm{Var}(X^{\bar{\pi}^{}}(\tau))$ is increasing with the parameter $\alpha>0$ for a given terminal time $\tau$ in Figure \ref{fig1}.

\begin{figure}[H]
\noindent\caption{Compare the efficient frontier with classical mean-variance model}
 \label{fig2}
\begin{center}
\includegraphics[width=4.2 in]{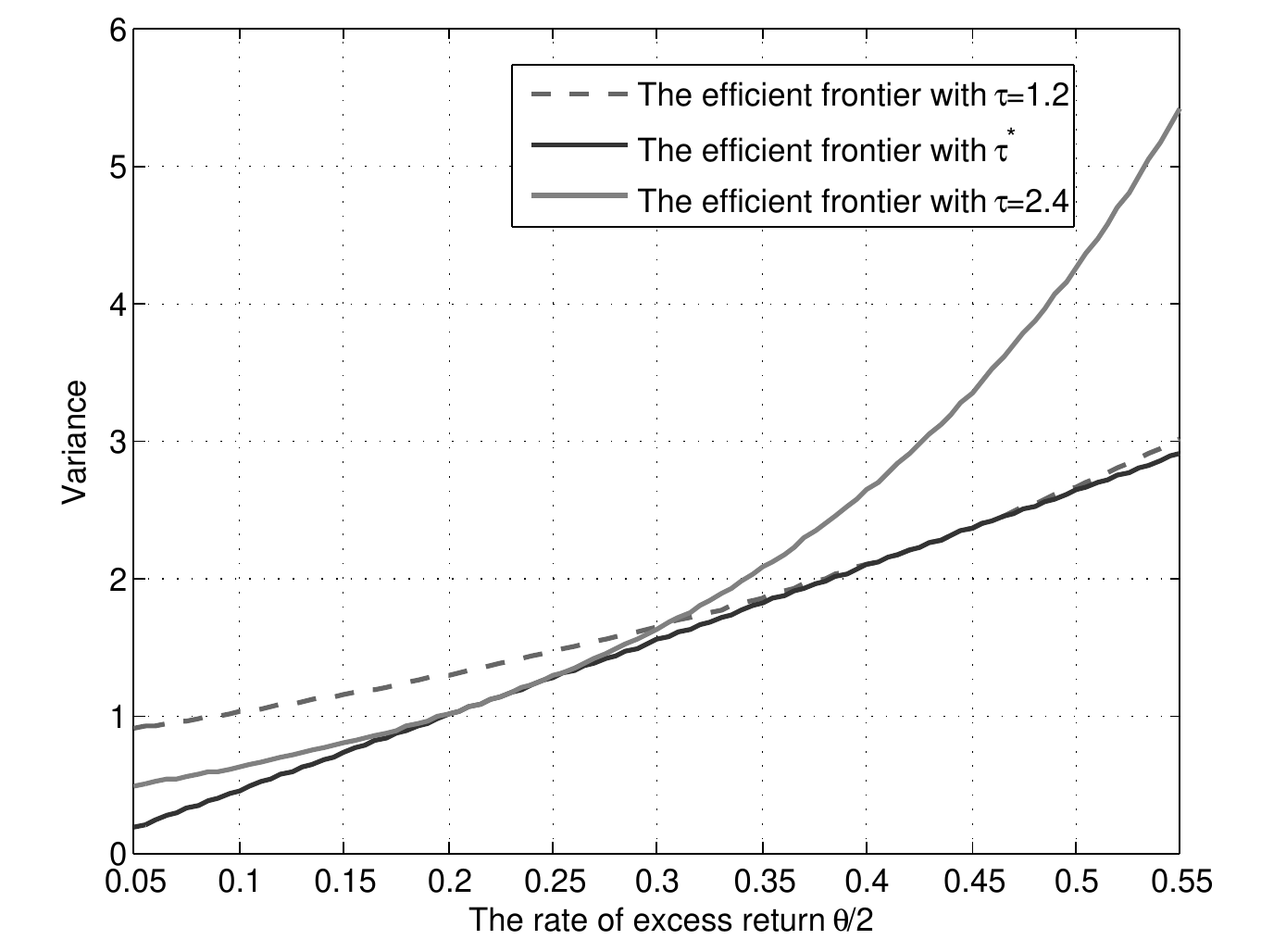}
\end{center}
\end{figure}
For the given parameters $r=5\%,\ b_1=10\%,\ \sigma_1=20\%$, we plot the relation between the rate of excess return and variance  as the efficient frontier. In Figure \ref{fig2}, for the given terminal time $\tau=1.2$, the related line  shows the relation of the rate of excess return $\displaystyle \frac{\theta}{2}\in [0.05,0.55]$ and the variance $\mathrm{Var}(X^{\bar{\pi}^{}}(\tau))$ of classical mean-variance model. For the  given terminal time $\tau=2.4$, the related line  shows the relation of the rate of excess return $\displaystyle \frac{\theta}{2}\in [0.05,0.55]$ and the variance $\mathrm{Var}(X^{\bar{\pi}^{}}(\tau))$ of classical mean-variance model. For the $\tau^*$, the related line shows the relation of the rate of excess return $\displaystyle \frac{\theta}{2}\in [0.05,0.55]$ and the variance $\mathrm{Var}(X^{{\pi}^{*}}(\tau^*))$ of our varying terminal time mean-variance model, where $\tau^*=\displaystyle \frac{1}{\phi}\ln\frac{{\theta}}{{\theta}-\phi}$ varies with $\theta$. The figure shows that the line of $\mathrm{Var}(X^{{\pi}^{*}}(\tau^*))$  is always under the line with $\tau=2.4$ and  line with $\tau=1.2$. Therefore, the proposed  model (\ref{cost-1}) can help to determine an optimal terminal time $\tau^*$ that varies  with  $\theta$. In addition, for  $\tau=1.2$ and $2.4$, we find that the line with $\tau=1.2$ touches the  line at $\displaystyle \frac{\theta}{2}=0.43$ and the  line with $\tau=2.4$ touches the red line at $\displaystyle \frac{\theta}{2}=0.225$, respectively. Therefore,  the variance $\mathrm{Var}(X^{\bar{\pi}^{}}(\tau))$ of classical mean-variance model takes the minimum  at some rate of excess return $\displaystyle \frac{\theta}{2}$  for the given terminal time $\tau$.

\begin{figure}[H]
 \caption{The optimal terminal time $\tau^*$ along with  the rate of excess return}
 \label{fig3}
\begin{center}
\includegraphics[width=4.2 in]{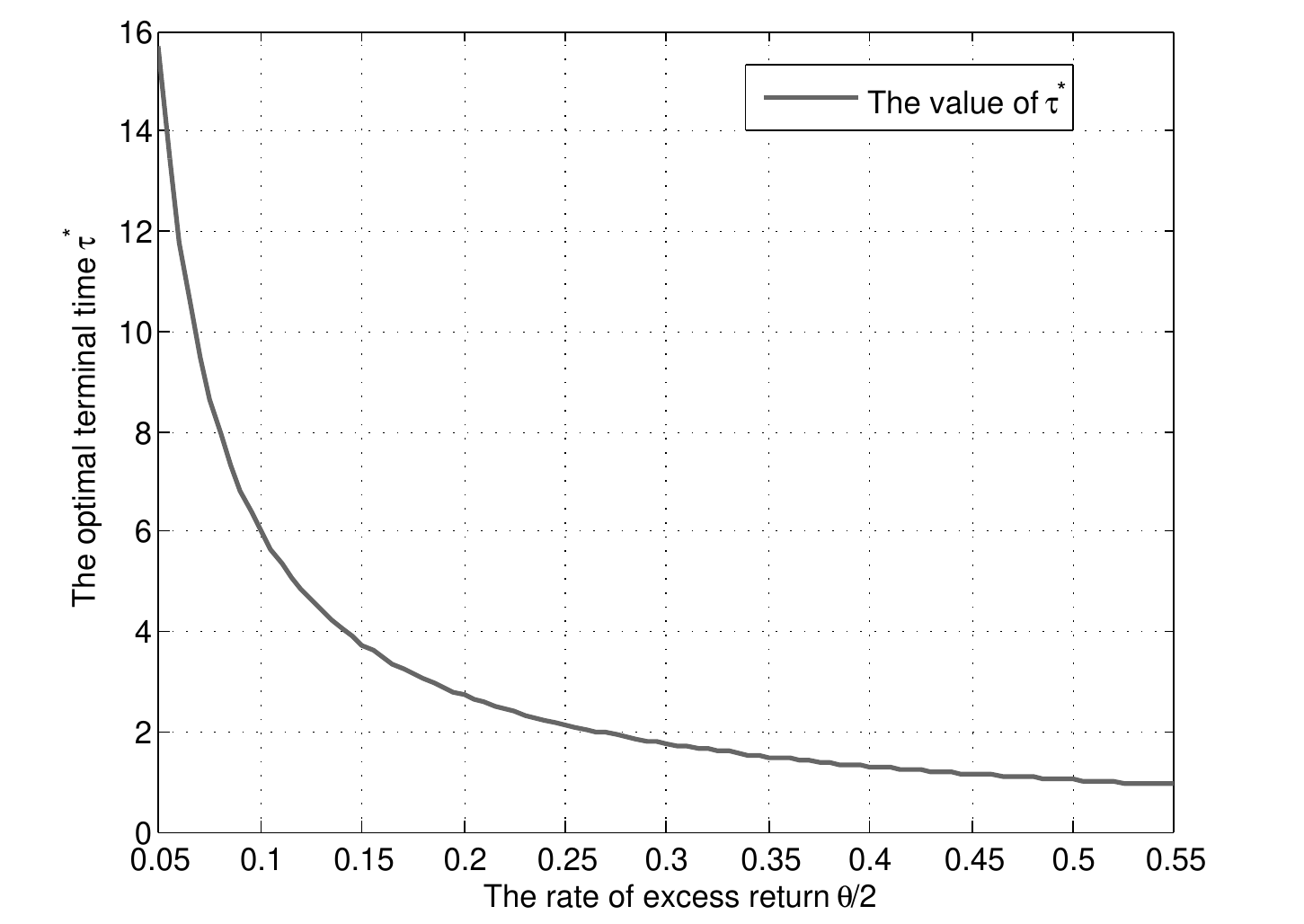}
\end{center}
\end{figure}
For the given parameters $r=5\%,\ b_1=10\%,\ \sigma_1=20\%$, we plot the relation between $\tau^*$ and the rate of  excess return $\displaystyle \frac{\theta}{2}\in [0.05,0.55]$ in Figure \ref{fig3}.
Here, $\tau^*=\displaystyle \frac{1}{\phi}\ln\frac{{\theta}}{{\theta}-\phi}$  is decreasing with $\displaystyle \frac{\theta}{2}\in [0.05,0.55]$. In Figure \ref{fig2}, the variance
$\mathrm{Var}(X^{{\pi}^{*}}(\tau^*))$  is increasing with $\displaystyle \frac{\theta}{2}\in [0.05,0.55]$. These results show that if we consider a small rate of excess return in the investment plan, for example $\displaystyle \frac{\theta}{2}=0.05$, then, we can keep the optimal strategy $\pi^*$ until $\tau^*=15.69$ with the variance $\mathrm{Var}(X^{{\pi}^{*}}(\tau^*))=0.72$. However, if we consider a high rate of excess return in the investment plan, for example $\displaystyle \frac{\theta}{2}=0.55$, then, we can keep the optimal strategy $\pi^*(\cdot)$ until $\tau^*= 0.94$ with the variance $\mathrm{Var}(X^{{\pi}^{*}}(\tau^*))=11.62$.

\bigskip

{ Now, we consider the relation between  the number of  risky assets $n$ and  $\tau^*$, let
$$
\tau^*=\displaystyle \frac{1}{\phi}\ln\frac{{\theta}}{{\theta}-\phi},\quad  \phi=\sum_{i=1}^n\bigg{(}\frac{b_i-r}{\sigma_i}\bigg{)}^2.
$$
\begin{figure}[H]
 \caption{The optima terminal time $\tau^*$ along with the number of risky asset}
 \label{fig4}
\begin{center}
\includegraphics[width=4.2 in]{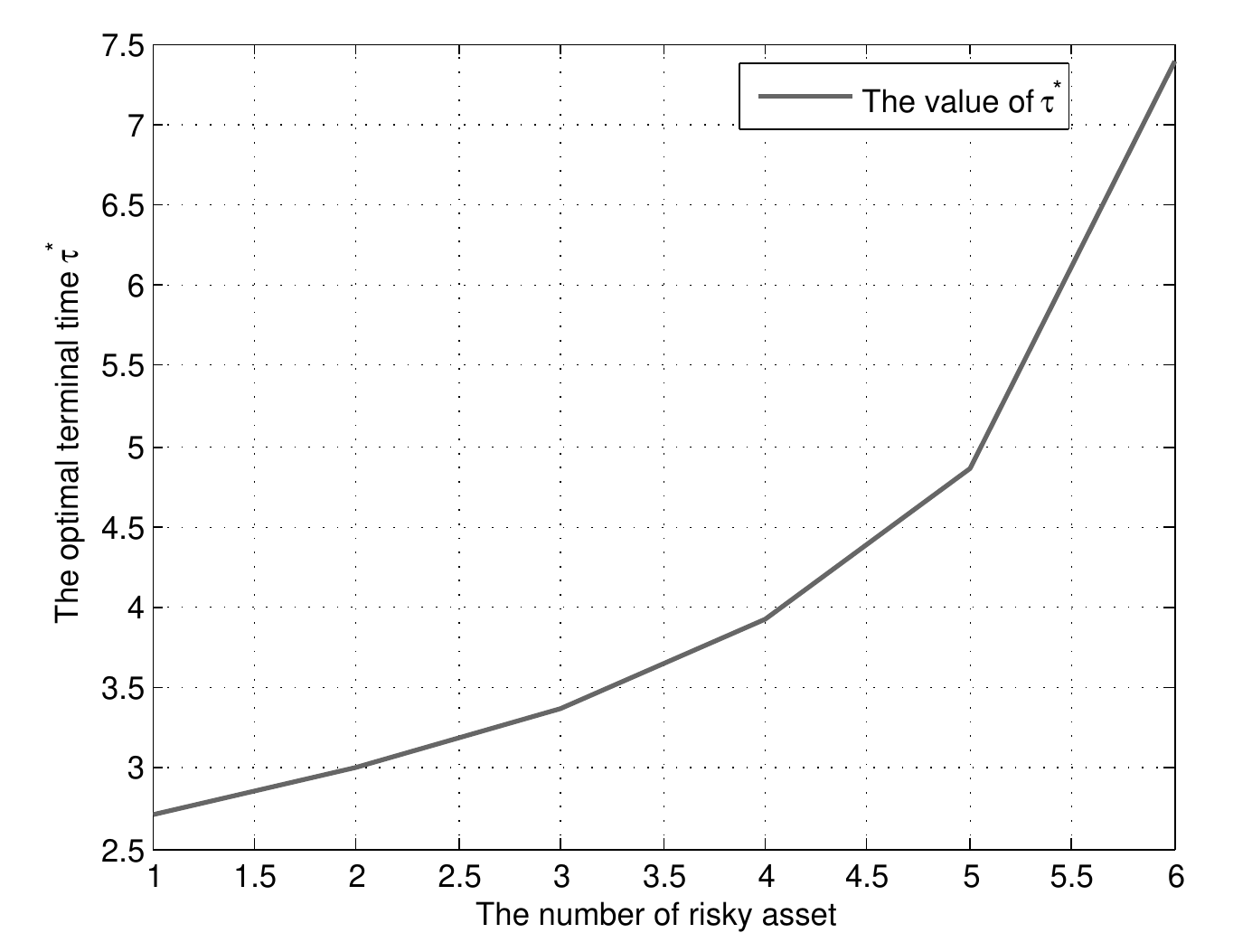}
\end{center}
\end{figure}
For the given parameters $r=5\%,\ b_i=10\%,\ \sigma_i=20\%,\ 1\leq i\leq 6,\ \theta=40\%$, we plot the relation between $\tau^*$ and $n$ in Figure \ref{fig4}.
Let $\theta>\phi$, $\tau^*=\displaystyle \frac{1}{\phi}\ln\frac{{\theta}}{{\theta}-\phi}$  is increasing with the number of risky assets, which indicates that the risk (variance) is decreasing with the number of risky assets. If $\theta\leq \phi$, Lemma \ref{le-5} indicates that we cannot obtain an optimal terminal time.

}

\section{Conclusion}
To improve the efficient frontier of the mean-variance model, we propose a varying terminal time mean-variance model with a constraint on  the varying terminal time. In the proposed model, we suppose that the investor's target moves with  the rate of return.

Our main results are as follows:
\begin{itemize}
\item  We  minimize the variance of the assets in a  portfolio, thus incorporating   the advantages of the classical mean-variance model.

 \item  The constraint on the varying terminal time allows us to find the optimal deterministic terminal time and to minimize the variance.

 \item The results of Section 4 show that the optimal terminal time is decreasing with the rate of excess return which suggests that the investor should change the holding time of an asset according to the rate of excess return.

 \item The proposed varying terminal time  mean-variance model  improves the efficient frontier of the classical mean-variance model.

\end{itemize}

This study represents the first step in considering the varying terminal time mean-variance problem, based on which, we can further investigate  topics  such as bankruptcy prohibition, transaction costs, and random parameters in complete and incomplete markets.


\appendix

\section{Proofs:}

\noindent \textbf{Proof of Lemma \ref{le-3}:} For any given $\tau> 0$, recall that
$$
\displaystyle \mathbb{E}[{X}^{\bar{\pi}^{}}(t)]=
xe^{\int_0^t(r(s)-\phi(s))
\mathrm{d}s}+\gamma(\tau)e^{-\int_t^{\tau}r(s)
\mathrm{d}s}[1-e^{-\int_0^{t}\phi(s)\mathrm{d}s}],\ t\leq \tau.
$$
The derivative of $\mathbb{E}[{X}^{\bar{\pi}^{}}(t)]$ at $t$ is given as follows,
\begin{equation*}
\begin{array}
[c]{rl}%
\displaystyle \frac{\mathrm{d}\mathbb{E}[{X}^{\bar{\pi}^{}}(t)]}{\mathrm{d}t}
&=\bigg{(}  xe^{\int_0^t(r(s)-\phi(s))
\mathrm{d}s}+\gamma(\tau)e^{-\int_t^{\tau}r(s)
\mathrm{d}s}[1-e^{-\int_0^{t}\phi(s)\mathrm{d}s}]\bigg{)}'\\
&>e^{-\int_0^{t}\phi(s)\mathrm{d}s}\bigg{[}x (r(t)-\phi(t))e^{\int_0^tr(s)
\mathrm{d}s} +\gamma(\tau)e^{-\int_t^{\tau}r(s)
\mathrm{d}s}\phi(t)\bigg{]}.
\end{array}
\end{equation*}
From
$$
\displaystyle \gamma(\tau)e^{-\int_t^{\tau}r(s)
\mathrm{d}s}=
\frac{1}{\mu(\tau)}e^{\int_0^{\tau}\phi(t)\mathrm{d}t-\int_t^{\tau}r(s)
\mathrm{d}s}+x e^{\int_0^{t}r(s)\mathrm{d}s},
$$
it follows that
$$
x (r(t)-\phi(t))e^{\int_0^tr(s)
\mathrm{d}s} +\gamma(\tau)e^{-\int_t^{\tau}r(s)
\mathrm{d}s}\phi(t)>0,
$$
and
$$
\displaystyle \frac{\mathrm{d}\mathbb{E}[{X}^{\bar{\pi}^{}}(t)]}{\mathrm{d}t}>0,\ t\leq \tau.
$$
Thus, $\mathbb{E}[{X}^{\bar{\pi}^{}}(\cdot)]$ is increasing on $[0,\tau]$, $\mathbb{E}[{X}^{\bar{\pi}^{}}(\tau)]=xh(\tau)$,  and we have
 $$
\tau={\tau}^{\bar{\pi}^{}}=\inf\bigg{\{}t:\mathbb{E}[{X}^{\bar{\pi}^{}}(t)]\geq xh(\tau),\ t\in [0,\tau] \bigg{\}}.
$$
Note that,
$$
\bar{\pi}^{}(t)=[\sigma(t)\sigma(t)^{\top}]^{-1}\beta(t)^{\top}(\gamma(\tau) e^{-\int_t^{\tau}r(s)\mathrm{d}s}-X^{\bar{\pi}^{}}(t))\in  L^2_{\mathcal{F}}[0,\tau;\mathbb{R}^n].
$$
Thus $\bar{\pi}^{}(\cdot)\in \mathcal{A}$ and $\mathcal{A}\neq \varnothing$, which completes this proof. $\ \ \ \ \ \ \ \ \Box$

\bigskip

\noindent \textbf{The Proof of Lemma \ref{le-5}:}
  By Theorem \ref{the-1}, for any given deterministic  time $\tau> 0$,  we have
\begin{equation}\label{apb-eq-1}
\displaystyle  \mathrm{Var}(X^{\bar{\pi}^{}}(\tau))=
\frac{x^2\bigg{(}h(\tau)
-e^{\int_0^{\tau}r(t)\mathrm{d}t}\bigg{)}^2}{e^{\int_0^{\tau}\phi(t)\mathrm{d}t}-1}.
\end{equation}
Based on  the decomposition of $xh(\tau)$ (\ref{eqh-0}), we set
$$
g(\tau)=(h(0)-1) e^{\int_0^{\tau}\frac{\theta(t)}{2}\mathrm{d}t}.
$$
 In the following, we calculate the derivative of $\mathrm{Var}(X^{\bar{\pi}^{}}(\tau))$ at $\tau$. Combining equations (\ref{apb-eq-1}) and $g(\tau)$, we have
\begin{equation*}
\begin{array}
[c]{rl}%
\displaystyle \frac{\mathrm{d}\mathrm{Var}(X^{\bar{\pi}^{}}(\tau))}{\mathrm{d}\tau}
=&\displaystyle \bigg{(} \displaystyle \frac{x^2\big{(}h(\tau)
-e^{\int_0^{\tau}r(t)\mathrm{d}t}\big{)}^2}{e^{\int_0^{\tau}\phi(t)\mathrm{d}t}-1}\bigg{)}'\\
=& \displaystyle \bigg{(}\frac{x(h(0)-1)e^{\int_0^{\tau}{\frac{{\theta}(t)}{2}}\mathrm{d}t}
}{e^{\int_0^{\tau}\phi(t)\mathrm{d}t}-1} \bigg{)}^2 \bigg{[}({\theta}(\tau)-\phi(\tau))e^{\int_0^{\tau}\phi(t)\mathrm{d}t} -{\theta}(\tau)\bigg{]}.\\
\end{array}
\end{equation*}
We set
$
I(\tau)=({\theta}(\tau)-\phi(\tau))e^{\int_0^{\tau}\phi(t)\mathrm{d}t} -{\theta}(\tau)
$
and have $I(0)=-\phi(0)<0$ from Assumption $\textbf{H}_2$.  By Assumption $\textbf{H}_3$, we have
 that $g(s)>0$ and
$$
{\theta}(s)=\frac{2g(s)'}{g(s)}=\displaystyle \frac{2h'(s)-2r(s)e^{\int_0^{s}r(t)\mathrm{d}t}}
{h(s)-e^{\int_0^{s}r(t)\mathrm{d}t}},\ s\geq 0.
$$

If ${\theta}(s)>\phi(s)+\delta, \ s\geq 0$ for a constant $\delta>0$, there exists $\tau_1$ and $\tau_2$ such that $I(\cdot)<0$ on $(0,\tau_1]$ and $I(\cdot)>0$ on $[\tau_2,+\infty)$. Thus,   $\mathrm{Var}(X^{\bar{\pi}^{}}(\cdot))$ is decreasing on $(0,\tau_1]$ and increasing on $[\tau_2,+\infty)$. By Assumption $\textbf{H}_1$,  we have that $I(\cdot)$ is a continuous function on $(0,+\infty)$. Thus there exists a minimum time $\tau^*\in[\tau_1,\tau_2]$, such that
$$
\displaystyle  \mathrm{Var}(X^{{\pi}^{*}}(\tau^*))=\min_{\tau\in(0,+\infty)}\mathrm{Var}(X^{\bar{\pi}^{}}(\tau)).
$$

If ${\theta}(s)\leq \phi(s), \ s\geq 0$, we have $I(\cdot)<0$ on $(0,+\infty)$ and
$$
\displaystyle  \lim_{\tau\to+\infty}\mathrm{Var}(X^{\bar{\pi}^{}}(\tau))=\displaystyle \inf_{\tau\in(0,+\infty)}
\mathrm{Var}(X^{\bar{\pi}^{}}(\tau)).
$$
This completes the proof. $\ \ \ \ \ \ \ \ \Box$

\bigskip

\noindent \textbf{The Proof of Theorem \ref{the-2}:}  We first consider the case $(i)$. Applying Lemma \ref{le-5}, we have that if
$$
\theta(s)> \phi(s)+\delta,\ s\geq 0,
$$
for a constant $\delta>0$, then there exists $(\pi^*(\cdot),\tau^*)$ that yields  the minimum value of $\mathrm{Var}(X^{\bar{\pi}^{}}(\tau))$ for any $\tau\in (0,+\infty)$. For any given $\tau> 0$, based on Theorem \ref{the-1}, it follows that $(\bar{\pi}^{}(\cdot),\tau)$ yields the minimum values of $\mathrm{Var}(X^{{\pi}}(\tau))$ for any $\pi(\cdot)\in \mathcal{A}$ with $\tau^{\pi}=\tau$. Note that, if $\tau^{\pi}<\tau$, we have $\mathbb{E}[X^{\tau^{\pi}}(\tau^{\pi})]=xh(\tau)$. For the given $\tau^{\pi}$, Theorem \ref{the-1} shows that there exists optimal strategy $\hat{\pi}(\cdot)$ such that
\begin{equation*}
\displaystyle  \mathrm{Var}(X^{\hat{\pi}}(\tau^{\pi}))=
\frac{\bigg{(}xh(\tau)
-xe^{\int_0^{\tau^{\pi}}r(t)\mathrm{d}t}\bigg{)}^2}{e^{\int_0^{\tau^{\pi}}\phi(t)\mathrm{d}t}-1}.
\end{equation*}
From Assumption $\textbf{H}_3$,
$
xh(\tau)>xe^{\int_0^{\tau^{}}r(t)\mathrm{d}t}>xe^{\int_0^{\tau^{\pi}}r(t)\mathrm{d}t},
$
we have
\begin{equation*}
\displaystyle  \mathrm{Var}(X^{\hat{\pi}}(\tau^{\pi}))
>
\frac{\bigg{(}xh(\tau^{})
-xe^{\int_0^{\tau^{}}r(t)\mathrm{d}t}\bigg{)}^2}{e^{\int_0^{\tau^{}}\phi(t)\mathrm{d}t}-1}
=\mathrm{Var}(X^{\bar{\pi}^{}}(\tau^{})).
\end{equation*}
Thus, $(\pi^*(\cdot),\tau^*)$ yields the minimum value of $\mathrm{Var}(X^{{\pi}}(\tau))$ for any $\tau\in (0,+\infty)$ and $\pi(\cdot)\in \mathcal{A}$. By Lemma \ref{le-3}, we have that
$\pi^*(\cdot)\in \mathcal{A}$, which shows that $(\pi^*(\cdot),\tau^*)$  solves  model (\ref{cost-1}) with  constraint (\ref{time-1}).

Next, we consider the case $(ii)$. If
$$
\theta(s)\leq \phi(s) ,\ s\geq 0.
$$
Then,  similarly to  case $(i)$, we have that $(\bar{\pi}^{}(\cdot),\tau)$ yields the minimum values of $\mathrm{Var}(X^{{\pi}}(\tau))$ for any $\pi(\cdot)\in L^2_{\mathcal{F}}[0,\tau;\mathbb{R}^n]$ with $\mathbb{E}[{X}^{{\pi}}(\tau)]=xh(\tau)$
and $\bar{\pi}^{}(\cdot)\in \mathcal{A}$. Thus, for a given $\tau>0$, $(\bar{\pi}^{}(\cdot),\tau)$ admits the minimum values of   model (\ref{cost-1}) for any $\pi(\cdot)\in \mathcal{A}$ with $\tau^{\pi}\leq \tau$ and cost functional (\ref{cost-1}) is decreasing with $\tau \in (0,+\infty)$. Thus, model (\ref{cost-1}) does not provide an optimal strategy in finite time. This completes the proof. $\ \ \ \ \ \ \ \ \Box$

\section*{Acknowledgement}
The author would like to thank the editor and two anonymous reviewers for their careful reading and suggestions for improving the quality of this paper.

\end{document}